\documentclass{article}

\usepackage{amsmath,amsfonts,amssymb}
\usepackage{graphicx}
\usepackage{caption}
\usepackage{subcaption}
\usepackage[top=2cm, left=2cm, bottom=2cm, right=2cm]{geometry}
\usepackage{epstopdf}
\usepackage{hyperref}
\usepackage{color}
\usepackage{cite}
\usepackage[utf8]{inputenc}

\begin{document}

\title{Discrete fractional Fourier transform: 
\\ Vandermonde approach}

\author{  H\'ector M. Moya-Cessa and Francisco Soto-Eguibar  \\
{\small Instituto Nacional de Astrof\'{\i}sica, \'Optica y Electr\'onica, INAOE} \\
\small {Calle Luis Enrique Erro 1, Santa Mar\'{\i}a Tonantzintla, Puebla, 72840 Mexico}\\
$^*$\small {Corresponding author: feguibar@inaoep.mx}}

\date{\today}

\maketitle

\begin{abstract}
Based on the definition of the Fourier transform in terms of the number operator of the quantum harmonic oscillator and in the corresponding definition of the fractional Fourier transform, we have obtained the discrete fractional Fourier transform from the discrete Fourier transform in a completely analogous manner. To achieve this, we have used a very simple method based on the Vandermonde matrices, to obtain rational and irrational powers of the discrete Fourier transform matrices.
\end{abstract}

\section{Introduction}
The discrete fractional Fourier transform (DFrFT) has been the subject of recent interest. Pei {\it et al.} proposed \cite{Pei} based on orthogonal projections. Since then a number of authors \cite{Candan,Singh,Sejdic,Serbes,Hanna} have studied the problem of defining the DFrFT by methods that involve mainly eigenvectors of the discrete Fourier transform matrix \cite{Mc,dickinson} because the additivity property of the DFrFT matrix requires having orthonormal eigenvectors of the Discrete Fourier transform $N\times N$ matrix, $F$. \\
The discrete Fourier transform may be implemented using circular waveguide arrays and may be related to a fundamental problem of quantum mechanics, namely, the search for a phase operator \cite{Leija}. On the DFrFT may be implemented also in such optical systems, but in this case, linear arrays with particular interaction constants \cite{Weimann}. The purpose of this contribution is to show that, by using a Vandermonde approach \cite{Champi} one can find functions of any given (square) matrix \cite{viana} and therefore can define the DFrFT matrix in an {\it easy} way.

\section{The discrete fractional Fourier transform}
In this work, we define the discrete fractional Fourier transform as
\begin{equation} \label{def}
F_\alpha \doteq U^\alpha
\end{equation}
where $\alpha$, a real number, is the order of the DFrFT and $U$ is the discrete Fourier transform matrix, given by
\begin{equation}\label{dft}
U_{j,k}=\frac{1}{\sqrt{N}}\left[\exp\left(-i\frac{2\pi}{N} \right)  \right]^{ \left( j-1\right) \left(k-1 \right)} ; \;
j=1,2,3,...,N, \;  k=1,2,3,...,N,
\end{equation}
being $N$ the dimension of the corresponding space. It is clear from the preceding equations that $F_\alpha F_{-\alpha}=1_{N\times N}$ with $1_{N\times N}$ the unit matrix, and that $F_1=U$.\\

To make definition \eqref{def} truly operational, we use the fact that given a complex vector space $V$ of finite dimension $N$ and any function $f$ well behaved, then $f(\mathbb{A})=\sum_{n=0}^{N-1} c_{n+1} \mathbb{A}^n$ for all $\mathbb{A} \in L\left( V \right)$, with the coefficients given by $\vec{c}=\left( {\mathbb{V}^{T}}\right) ^{-1} \vec{f} \left( \lambda \right)$, where $\mathbb{V} \in L\left( V \right)$ is the confluent Vandermonde matrix $N \times N$ associated with $\mathbb{A}$ \cite{viana}.
The confluent Vandermonde matrix associated with $\mathbb{A}$ is constructed as follows: Let $\left\lbrace \lambda_i, i=1,2,3,...,n  \right\rbrace $ the $n$ different eigenvalues of the matrix $\mathbb{A}$ each one of multiplicity $m_i$; note that $\sum_{i=1}^{n}m_i=N$. The confluent Vandermonde matrix associated with $\mathbb{A}$ is given by
\begin{equation}
\mathbb{V}=\left(
\begin{array}{ccccc}
\mathbb{V}_1 & \mathbb{V}_2 & \mathbb{V}_3&...& \mathbb{V}_n
\end{array}
\right),
\end{equation}
where $\mathbb{V}_i, \; i=1,2,3,...,n$, are matrices $N \times m_i$ given by
\begin{equation}
\left( \mathbb{V}_i\right)_{j,k}=
\left[ \frac{d^{k-1}\lambda^{j-1}}{d\lambda^{k-1}}\right]_{\lambda=\lambda_i}
\; j=1,2,3,...,N, \, k=1,2,3,...,m_i;  \; i=1,2,3,...,n.
\end{equation} 
The vector $\vec{f}\left( \lambda\right) $ stands for
\begin{equation}
\vec{f}\left( \lambda\right)=
\left( 
\begin{array}{c}
\vec{f}_1\\
\vec{f}_2\\
\vec{f}_3\\
\vdots\\
\vec{f}_n\\
\end{array}
\right) .
\end{equation}
where
\begin{equation}
\left( \vec{f}_i\right)_j=\left[ \frac{d^{j-1} f\left(\lambda \right) }{\lambda^{j-1}}\right]_{\lambda=\lambda_i};
\; j=1,...,m_i, \; i=1,2,3,...,n.
\end{equation}\\

In the case of the finite Fourier transform matrix, Eq. \eqref{dft}, the characteristic polynomial is \cite{Mc,dickinson}
\begin{equation}
p\left(\lambda \right)=\left(\lambda-1 \right)^{\lfloor\frac{N+4}{4}\rfloor} \left( \lambda+1\right)^{\lfloor\frac{N+2}{4}\rfloor}  \left(\lambda+i \right)^{\lfloor\frac{N+1}{4}\rfloor}  \left(\lambda-i\right)^{\lfloor\frac{N-1}{4}\rfloor}   ,
\end{equation} 
so the eigenvalues are the fourth roots of unity, ${+1,-1,-i,+i}$, and their multiplicities are presented in the following table \cite{Mc,dickinson}:
\begin{table}[h!]
\centering
\label{tabla}
\begin{tabular}{l||l|l|l|l}
	$N$&$\lambda=+1$&$\lambda=-1$&$\lambda=-i$&$\lambda=+i$\\
	\hline
	$4m$&$m+1$&$m$&$m$&$m-1$\\
	$4m+1$&$m+1$&$m$&$m$&$m$\\
	$4m+2$&$m+1$&$m+1$&$m$&$m$\\
	$4m+3$&$m+1$&$m+1$&$m+1$&$m$\\
\end{tabular}
\caption{Table of the multiplicity of the eigenvalues of the DFT matrix}
\end{table}\\
The confluent Vandermonde matrix associated with the finite Fourier transform matrix \eqref{dft} is then given by \cite{Champi, viana}
\begin{equation}
\Lambda=\left(
\begin{array}{cccc}
\Lambda_1 & \Lambda_2 & \Lambda_3& \Lambda_4
\end{array}
\right),
\end{equation}
where
\begin{equation}
\left( \vec{\Lambda}_i\right)_{k,j}=\left[ \frac{d^{j-1}\lambda^{k-1}}{\lambda^{j-1}}\right]_{\lambda=\lambda_i}
\; j=1,...,m_i, \, k=1,2,...,N;  \; i=1,2,3,4.
\end{equation}
The integer numbers $\left\lbrace m_1,m_2,m_3,m_4\right\rbrace $ are the corresponding multiplicities of the fourth different possibles eigenvalues $\left\lbrace +1,-1,-i,+i\right\rbrace $ given in Table 1.\\
The vector $\vec{f}\left( \lambda\right) $ is given in this particular case by
\begin{equation}
\vec{f}\left( \lambda\right)=
\left( 
\begin{array}{c}
\vec{f}_1\\
\vec{f}_2\\
\vec{f}_3\\
\vec{f}_4\\
\end{array}
\right) .
\end{equation}
where
\begin{equation}
\left( \vec{f}_i\right)_j=\left[ \frac{d^{j-1}\lambda^\alpha}{\lambda^{j-1}}\right]_{\lambda=\lambda_i}
\; j=1,...,m_i, \; i=1,2,3,4.
\end{equation}
It is well known that the inverse of the transpose of the confluent Vandermonde matrix always exists  and explicit expressions have been found for it \cite{Champi}.

\section{An example}
As a practical example, we will construct in this section $F_{3/7}$ with $N=7$. The corresponding confluent Vandermonde matrix is in this case
\begin{equation}
\Lambda=
\left(
\begin{array}{ccccccc}
 1 & 0 & 1 & 0 & 1 & 0 & 1 \\
 1 & 1 & -1 & 1 & -i & 1 & i \\
 1 & 2 & 1 & -2 & -1 & -2 i & -1 \\
 1 & 3 & -1 & 3 & i & -3 & -i \\
 1 & 4 & 1 & -4 & 1 & 4 i & 1 \\
 1 & 5 & -1 & 5 & -i & 5 & i \\
 1 & 6 & 1 & -6 & -1 & -6 i & -1 \\
\end{array}
\right).
\end{equation}
The inverse of the transpose confluent Vandermonde matrix results to be
\begin{equation}
\left( \Lambda^T\right)^{-1}= 
\left(
\begin{array}{ccccccc}
 \frac{1}{4}+\frac{3 i}{16} &
   -\frac{1}{16}-\frac{i}{16} &
   \frac{1}{4}-\frac{3 i}{16} &
   \frac{1}{16}-\frac{i}{16} & \frac{7}{16} &
   \frac{i}{8} & \frac{1}{16} \\
 \frac{1}{2} & -\frac{1}{8} & -\frac{1}{2} &
   -\frac{1}{8} & \frac{i}{8} & 0 & -\frac{i}{8}
   \\
 \frac{3}{8}-\frac{3 i}{16} &
   -\frac{1}{16}+\frac{i}{16} &
   \frac{3}{8}+\frac{3 i}{16} &
   \frac{1}{16}+\frac{i}{16} & -\frac{9}{16} &
   -\frac{i}{8} & -\frac{3}{16} \\
 \frac{1}{4} & 0 & -\frac{1}{4} & 0 &
   -\frac{i}{4} & 0 & \frac{i}{4} \\
 -\frac{3 i}{16} & \frac{1}{16}+\frac{i}{16} &
   \frac{3 i}{16} & -\frac{1}{16}+\frac{i}{16} &
   -\frac{3}{16} & -\frac{i}{8} & \frac{3}{16}
   \\
 -\frac{1}{4} & \frac{1}{8} & \frac{1}{4} &
   \frac{1}{8} & \frac{i}{8} & 0 & -\frac{i}{8}
   \\
 -\frac{1}{8}+\frac{3 i}{16} &
   \frac{1}{16}-\frac{i}{16} &
   -\frac{1}{8}-\frac{3 i}{16} &
   -\frac{1}{16}-\frac{i}{16} & \frac{5}{16} &
   \frac{i}{8} & -\frac{1}{16} \\
\end{array}
\right)
\end{equation}
The $\vec{f}$ vector is
\begin{equation}
\vec{f}=
\left(1,\frac{3}{7},(-1)^{3/7},-\frac{3}{7}
   (-1)^{3/7},-(-1)^{11/14},\frac{3}{7}
   (-1)^{2/7},(-1)^{3/14}\right).
\end{equation}
It is possible now to calculate $\vec{c}=\left( {\mathbb{V}^{T}}\right) ^{-1} \vec{f} \left( \lambda \right)$, and finally  $F_{3/7}=\sum_{n=0}^{6} c_{n+1} U^n$. The matrix obtained is too big to be displayed here. However, as a verification that the matrix corresponding to the discrete fractional Fourier transform behaves correctly,  we checked that $F_{3/7}^{7/3}$ reproduces the discrete Fourier matrix, $U$. Also, if we calculated $F_{3/7}^{28/3}$  obtaining the identity matrix, as expected.

\section{Conclusions}
Based on the fact that the Fourier transform may be written as
\begin{equation}
{\mathcal F}=e^{i\frac{\pi}{2}a^{\dagger}a},
\end{equation}
with $a$ and $a^{\dagger}$ the annihilation and creation operators \cite{Soto} for the harmonic oscillator (see for instance \cite{Fan}) and that
the fractional Fourier transform of $s$ order is obtained directly from it
\begin{equation}
{\mathcal F}_s= \left(e^{i\frac{\pi}{2}a^{\dagger}a}\right)^s,
\end{equation}
we have obtained the discrete fractional Fourier transform from the discrete Fourier transform in an analogous manner.
We have used a simple method based on Vandermonde matrices \cite{viana} to obtain functions of a square matrix when it has (some) equal eigenvalues, as is the case for the discrete Fourier transform.

\end{document}